\documentclass[10pt,reqno]{amsart}

\usepackage{amssymb}
\usepackage{amsthm}
\usepackage{amsmath}

\newtheorem{thm}{Theorem}
\newtheorem{cor}[thm]{Corollary}
\newtheorem{lem}[thm]{Lemma}
\newtheorem{knownthm}{Theorem}

\theoremstyle{definition}
\newtheorem{rem}{Remark}[section]

\newcommand{\round}{\partial}
\newcommand{\CC}{\widehat{\mathbb{C}}}
\newcommand{\C}{\mathbb{C}}
\newcommand{\D}{\mathbb{D}}
\newcommand{\R}{\mathbb{R}}

\newcommand{\A}{\mathcal{A}}

\newcommand{\pf}{\noindent{\bf\textit{Proof. \,}}}
\newcommand{\epf}{\hspace*{\fill}$\Box$}
\newcommand{\no}{\noindent}
\newcommand{\dstyle}{\displaystyle}
\renewcommand{\Re}{\textup{Re}\,}

\title[]{Ruscheweyh's univalence criterion and quasiconformal extensions}
\author[]{IKKEI HOTTA}
\address{Division of Mathematics, Graduate School of Information Sciences, Tohoku University, 6-3-09 Aramaki-Aza-Aoba, Aoba-ku, Sendai 980-8579, Japan}
\email{ikkeihotta@ims.is.tohoku.ac.jp}

\begin{document}
\maketitle

\small
\textbf{Abstract.}
Ruscheweyh extended the work of Becker and Ahlfors on sufficient conditions for a normalized analytic function on the unit disk to be univalent there.
In this paper we refine the result to a quasiconformal extension criterion with the help of Becker's method. 
As an application, a positive answer is given to an open problem proposed by Ruscheweyh.

\normalsize



\section{Introduction}
Throughout the paper, $\D$ denotes the unit disk $\{|z|<1\}$ in the complex plane $\C$ and $\D^{*}$ the exterior domain of $\D$ in the Riemann sphere $\CC = \C \cup \{\infty\}$.

Let $\A$ be a family of normalized analytic functions $f(z) = z + \sum_{n=2}^{\infty}a_{n}z^{n}$ on $\D$.
We say that a sense-preserving homeomorphism $f$ of a plane domain $G \subset \C$ is $k$-quasiconformal if $f$ is absolutely continuous on almost all lines parallel to the coordinate axes and $|f_{\bar{z}}| \leq k|f_{z}|$, almost everywhere $G$, where $f_{\bar{z}} = \round f/\round \bar{z}$, $f_{z} = \round f/\round z$ and $k$ is a constant with $0 \leq k <1$.

Ahlfors \cite{Ahlfors:1974} has shown that the following condition is sufficient for quasiconformal extensibility of univalent functions as an extension of Becker's univalence condition \cite{Becker:1972} (see also \cite{Pom:1975}, p175);

\begin{knownthm}[\cite{Ahlfors:1974},\cite{Becker:1976}]\label{AB}
Let $f \in \A$.
If there exists a $k,\,0 \leq k < 1,$ such that for a constant $c \in \C$ satisfying $|c| \leq k$ and all $z \in \D$
\begin{equation}\label{AB}
\left|
c |z|^{2} + (1-|z|^{2})\frac{zf''(z)}{f'(z)}
\right|
\leq k
\end{equation}
then $f$ has a $k$-quasiconformal extension to $\C$.
\end{knownthm}

The limiting case $k \to 1$ in the above theorem ensures univalence of $f$ in $\D$.
Ruscheweyh \cite{Ruscheweyh:1976} extended this univalence condition in the following way;

\begin{knownthm}[\cite{Ruscheweyh:1976}]\label{Ru}
Let $s = a+ib,\,a>0,\,b \in \R$ and $f \in \A$.
Assume that for a constant $c \in \C$ and all $z \in \D$
\begin{equation}\label{ruscheweyh}
\left| c|z|^{2}+s-a(1-|z|^{2})\left\{s\left(1+ \frac{zf''(z)}{f'(z)}\right)+(1-s)\frac{zf'(z)}{f(z)}\right\}\right|
\leq
M
\end{equation}
with
$$
M =
\left\{
\begin{array}{ll}
\dstyle a|s| + (a-1)|s+c|, &\textit{if}\hspace{10pt}   0 < a \leq 1, \\[5pt]
\dstyle |s|, &\textit{if}\hspace{10pt}   1 < a,
\end{array} 
\right.
$$
then $f$ is univalent in $\D$.
\end{knownthm}

The case $s=1$ with $c$ replaced by $-1-c$ is the special case of Theorem \ref{AB}.

The purpose of this paper is to refine Ruscheweyh's univalence condition to a quasiconformal extension criterion
which includes Theorem \ref{AB};

\begin{thm}\label{main}
Let $s = a+ib,\,a>0,\,b \in \R,\,k \in [0,1)$ and $f \in \A$.
Assume that for a constant $c \in \C$ and all $z \in \D$
\begin{equation}\label{maineq}
\left| 
c|z|^{2}+s-a(1-|z|^{2})\left\{s\left(1+ \frac{zf''(z)}{f'(z)}\right)+(1-s)\frac{zf'(z)}{f(z)}\right\}
\right|
\leq M
\end{equation}
with
$$
M =
\left\{
\begin{array}{ll}
a k |s| + (a-1)|s+c|, &\textit{if}\hspace{10pt}  0 < a \leq 1, \\[5pt]
k |s|, &\textit{if}\hspace{10pt}1 < a,
\end{array} 
\right.
$$
then $f$ has an $l$-quasiconformal extension to $\C$, where
\begin{equation}\label{mainl}
l=
\frac
{
2ka+(1-k^{2}) |b|
}
{
(1+k^{2})a+(1-k^{2}) |s|
}
<1.
\end{equation}
\end{thm}

\begin{rem}\label{rem11}
If $f \in \A$, then it is easy to verify that there exists a sequence $\{z_{n}\} \subset \D$ with $|z_{n}| \to 1$ such that for each $s \in \{z \in \C : \Re z >0\}$
$$
\sup_{n}\left|s\left(1+ \frac{z_{n}f''(z_{n})}{f'(z_{n})}\right)+(1-s)\frac{z_{n}f'(z_{n})}{f(z_{n})}\right| < \infty
$$
which shows that \eqref{maineq} implies the inequality
\begin{equation}\label{inequality}
|c+s|\leq M.
\end{equation}
This inequality is needed for proving that $f(z)$ has no zeros in $0<|z|<1$ (see Lemma \ref{approx}). 
In \cite{Ruscheweyh:1976}, it is mentioned that \eqref{maineq} implies $f(z) \neq 0$,\, $0<|z|<1$, without proof.
The part of \eqref{inequality} can be found in \cite{Ruscheweyh:1976}.

\begin{rem}
A similar argument to Remark \ref{rem11} is also valid for Theorem \ref{AB}.
It follows that the assumption $|c| \leq k$ is embedded in the inequality \eqref{AB}.
\end{rem}



\end{rem}

The next application follows from Theorem \ref{main}.
Let  $\alpha>0$ and $\beta \in \R$.
It follows from a result of Sheil-Small \cite[Theorem 2]{Sheil-Small:1972} that 
\begin{equation}\label{bazilevic}
\mathrm{Re} \left\{1+\frac{zf''(z)}{f'(z)} +(\alpha +i\beta-1)\frac{zf'(z)}{f(z)}\right\} >0 \hspace{20pt}(z \in \D)
\end{equation}
is sufficient for $f \in \A$ to be a Bazilevi\v c function of type ($\alpha, \beta$) \footnote{The author would like to thank Professor Yong Chan Kim for this remark.} (see also \cite{KimSugawa}).
Here, a function $f \in \A$ is called \textit{Bazilevi\v c of type ($\alpha, \beta$)} if
$$
f(z)=
\left[
(\alpha +i\beta)
\int_{0}^{z}
 g(\zeta)^{\alpha} h(\zeta) \zeta^{i\beta-1}d\zeta
\right]^{1/(\alpha+i \beta)}
$$
for a starlike univalent function $g \in \A$ and an analytic function $h$ with $h(0)=1$ satisfying $\textit{Re} (e^{i\lambda}h) > 0$ in $\D$ for some $\lambda \in \R$.
Together with this fact, the next theorem follows;

\begin{thm}\label{maincor}
Let $\alpha>0, \beta \in \R$ and $k \in [0,1)$.
If $f \in \A$ satisfies
\begin{equation}\label{maincorineq}
\left|
1+\frac{zf''(z)}{f'(z)}+(\alpha+i \beta-1)  \frac{zf'(z)}{f(z)}
-\frac{\alpha^{2}+\beta^{2}}{\alpha}
\right|
\leq M
\end{equation}
for all $z \in \D$ with
$$
M =
\left\{
\begin{array}{ll}
 k & \textit{if}\hspace{10pt}   \alpha < \alpha^{2}+\beta^{2}, \\[5pt]
k (\alpha^{2}+\beta^{2})/\alpha &\textit{if}\hspace{10pt}   \alpha^{2}+\beta^{2} \leq \alpha,
\end{array} 
\right.
$$
then $f$ is a Bazilevi\v c function of type ($\alpha, \beta$) and can be extended to a $\tilde{k}$-quasiconformal automorphism of $\C$, where
$$
\tilde{k} = 
\frac
{
2k\alpha+(1-k^{2}) |\beta|
}
{
(1+k^{2})\alpha+(1-k^{2})\sqrt{\alpha^{2}+\beta^{2}}
}.
$$
\end{thm}

Next, we shall discuss quasiconformal extensibility of functions $\dstyle g(z)= z+\frac{d}{z} + \cdots$ analytic in $\D^{*}$.



\begin{thm}\label{mainD*}
Let $s=a+ib,\,a \geq 1, b \in \R$ and $k \in [0,1)$ which satisfies $|b/s| \leq k$.
Let $\dstyle g(\zeta)= \zeta +\frac{d}{\zeta} + \cdots$ be analytic in $\D^{*}$ and fulfill
\begin{equation}\label{maineqD*}
\left|
ib+(1-|\zeta|^{2})a
\left
\{(1-s)
\left(
1-\frac{\zeta g'(\zeta)}{g(\zeta)}
\right)
-s \frac{\zeta g''(\zeta)}{g'(\zeta)}
\right\}
\right|
\leq
ak|s|-|b|(a-1)
\end{equation}
for all $\zeta \in \D^{*}$. 
Then $g$ can be extended to an $l$-quasiconformal automorphism of $\CC$, where
$$
l=
\frac
{
2ka+(1-k^{2}) |b|
}
{
(1+k^{2})a+(1-k^{2}) |s|
}.
$$
\end{thm}

The case $k \to 1$ corresponds to a univalence criterion which is due to Ruscheweyh \cite{Ruscheweyh:1976}.

Theorem \ref{mainD*} yields the following corollary which gives a positive answer to an open problem proposed by Ruscheweyh \cite{Ruscheweyh:1976}, i.e., whether a function $g(\zeta) = \zeta + d/\zeta + \cdots$ with $(|\zeta|^{2}-1)|1+ (\zeta f''(\zeta)/f'(\zeta))-(\zeta f'(\zeta)/f(\zeta)) | \leq k$ for all $\zeta \in \D^{*}$ admits a quasiconformal extension to $\C$;

\begin{cor}\label{corD*}
Let $\dstyle g(\zeta)= \zeta+\frac{d}{\zeta} + \cdots$ be analytic in $\D^{*}$.
If there exists $k \in [0,1)$ such that
$$
(|\zeta|^{2}-1)
\left|
1+\frac{\zeta g''(\zeta)}{g'(\zeta)} -\frac{\zeta g'(\zeta)}{g(\zeta)}
\right|
\leq k
$$
for all $\zeta \in \D^{*}$, then $g$ can be extended to a $k$-quasiconformal automorphism of $\CC -\{0\}$.
\end{cor}

From the above corollary we have another extension criterion for analytic functions on $\D$;

\begin{cor}\label{corD*D}
Let $f \in \A$ with $f''(0)=0$.
If there exists $k \in [0,1)$ such that
$$
(1-|z|^{2})
\left|
1+\frac{zf''(z)}{f'(z)} -\frac{zf'(z)}{f(z)}
\right|
\leq k
$$
for all $z \in \D$, then $f$ can be extended to a $k$-quasiconformal automorphism of $\C$.
\end{cor}


\

\section{Preliminaries}
Our investigations are based on the theory of L\"owner chains.
A function $f_{t}(z) =f(z,t) = a_{1}(t)z + \sum_{n=2}^{\infty}a_{n}(t)z^{n},\,a_{1}(t) \neq 0,$ defined on $\D \times [0,\infty)$ is called a \textit{L\"owner chain} if $f_{t}(z)$ is holomorphic and univalent in $\D$ for each $t \in [0,\infty)$ and satisfies $f_{s}(\D) \subsetneq f_{t}(\D)$ and $f(0,s) = f(0,t)$ for $0 \leq s \leq t < \infty$, and if $a_{1}(t)$ is locally absolutely continuous in $t \in [0,\infty)$ with $\lim_{t \to \infty}|a_{1}(t)| = \infty$.
Then $f(z,t)$ is absolutely continuous 
in $t \in [0,\infty)$ for each $z \in \D$ 
and satisfies the \textit{L\"owner differential equation}
\begin{equation}\label{loewner}
\dot{f}(z,t) =  h(z,t)z f'(z,t)
\end{equation}
for $z \in \D$ and almost every $t \in [0,\infty)$.
Here, $\dot{f}(z,t) = \round f(z,t) /\round t,\,f'(z,t)=\round f(z,t)/\round z$ and $h(z,t)$ is a function measurable on $t \in [0,\infty)$, holomorphic in $|z| < 1$ and $\mathrm{Re} h(z,t) > 0$ (\cite{Pom:1965}).



An interesting method connecting the theory of quasiconformal extensions with L\"owner chains was obtained by Becker;


\begin{knownthm}[\cite{Becker:1972}, see also \cite{Becker:1980}]\label{Beckerthm}
Suppose that $f(z,t)$ is a L\"owner chain for which $h(z,t)$ of \eqref{loewner} satisfies the condition 
$$
\left|
\frac{h(z,t)-1}{h(z,t)+1}
\right| \leq k
$$
Then $f_{t}(z)$ admits a continuous extension to $\overline{\D}$ for each $t \geq 0$ and the map defined by
$$
\hat{f}(re^{i\theta}) =
\left\{
\begin{array}{ll}
\dstyle f(re^{i\theta},0) & \textit{if}\hspace{10pt}  r < 1, \\[5pt]
\dstyle f(e^{i\theta},\log r) & \textit{if}\hspace{10pt}  r \geq 1,
\end{array} 
\right.
$$
is a $k$-quasiconformal extension of $f_{0}$ to $\C$.
\end{knownthm}

\

\section{Proof of Theorem \ref{main}}

The proof is divided into two parts. 
The first part of the proof is based on \cite{Ruscheweyh:1976}.

(i)
First we assume that $f(z)/z \neq 0$ for all $z \in \D$.
Then we can define
$$
f(z,t) = 
f(e^{-s t}z)
\left\{
1-\frac{a}{c}(e^{2t}-1)
\frac{e^{-s t}z f'(e^{-s t}z)}{f(e^{-s t}z)}
\right\}^{s}
$$
and let
\begin{equation}\label{chain2}
F(z,t) = f(z,t/|s|).
\end{equation}
A straightforward calculation shows
\begin{equation}\label{h(z,t)}
h(z,t)
= 
\frac{\dot{F}(z,t)}{zF'(z,t)}
=
\frac{s}{|s|}\cdot
\frac{1+P(e^{-st/|s|}z,t/|s|)}{1-P(e^{-st/|s|}z,t/|s|)},
\end{equation}
where
$$
P(z,t) = 
\frac{c}{a}e^{-2t} +1+(e^{-2t}-1)H_{s}(z)
$$
and
$$
H_{s}(z)=s\left(1+\frac{zf''(z)}{f'(z)}\right)    +    (1-s)\frac{zf'(z)}{f(z)}.
$$
Since $h(z,t)$ is holomorphic in $z \in \D$ and measurable on $t \in [0,\infty)$, applying Theorem \ref{Beckerthm} to \eqref{h(z,t)}, we see that the condition
$$
\left|
\frac
{s (1+P(e^{-st/|s|}z,t/|s|))    -    |s| (1-P(e^{-st/|s|}z,t/|s|))}
{s (1+P(e^{-st/|s|}z,t/|s|))    +    |s| (1-P(e^{-st/|s|}z,t/|s|))}
\right|
\leq l
$$
implies $l$-quasiconformal extensibility of $f(z)$.
This is equivalent to
\begin{equation}\label{diskA}
\left|
P+
\frac
{
(1+l^{2})b
}
{
(1+l^{2})a+(1-l^{2}) |s|
}
i
\right|
\leq
\frac
{
2l |s|
}
{
(1+l^{2})a+(1-l^{2}) |s|
}.
\end{equation}

Here, we shall prove the following Lemma;
\begin{lem}\label{lemdiskB}
Under the assumption of Theorem \ref{main}, we have 
\begin{equation}\label{diskB}
|a P(e^{-st/|s|}z,t/|s|) + ib| <  k|s|
\end{equation}
for $z \in \D$ and $t \in [0,\infty)$.
\end{lem}

\pf
We have
$$
|a P + ib| \leq m_{1}+m_{2}
$$
by triangle inequality, where
$$
m_{1}=
(1-e^{-2t/|s|})
\left|
\frac
{c e^{-2a t/|s|} + s}
{1-e^{-2a t/|s|}}
-a  H_{s}(e^{-st/|s|}z)
\right|
$$
and
$$
m_{2}=
\left|
(c e^{-2a t/|s|} +s)
\frac{1-e^{-2t/|s|}}{1-e^{-2a t/|s|}}
-(c e^{-2t/|s|} + s)
\right|.
$$
Then it is enough to show that $m_{1}+m_{2} < k|s|.$
(\ref{maineq}) implies
$$
\left|
\frac
{c   |e^{st/|s|}z|^{2}   +   s}
{1   -   |e^{st/|s|}z|^{2}}
-a  H_{s}(e^{-st/|s|}z)
\right|
\leq
\frac
{M}
{1   -   |e^{st/|s|}z|^{2}}
\leq
\frac
{M}
{1   -   e^{-2a t/|s|}}
$$
for $z \in \D$.
Let $q(t)=(1-e^{-2t/|s|})/(1-e^{-2a t/|s|})$. 
Applying the maximum modulus principle to the function
$$
\frac
{c e^{-2a t/|s|} + s}
{1-e^{-2a t/|s|}}
-a  H_{s}(e^{-st/|s|}z)
$$
we have $$m_{1} \leq q(t)M.$$
On the other hand $$m_{2} \leq |c+s||1-q(t)|.$$

\no
Since  $1 \leq q(t) < 1/a$ if $0 < a \leq 1$ and $1/a <q(t) \leq 1$ if $1 < a$ for all $t \in [0,\infty)$, we conclude that $m_{1}+m_{2} < k|s|$ which is our desired inequality.\epf

We now let $\Delta$ and $\Delta'$ be disks which are defined by replacing $P$ in \eqref{diskA} and \eqref{diskB} to a complex variable $w$.
 It remains to find the smallest $l$ so that $\Delta'$ is contained by $\Delta$. 
 Note that if $k=l=1$ then these two disks coincide.
The following condition is necessary and sufficient for $\Delta' \subset \Delta$;
\begin{equation}\label{Inequality01}
\left|
\frac
{
(1+l^{2})b
}
{
(1+l^{2})a+(1-l^{2}) |s|
}
-
\frac{b}{a}
\right|
\leq
\frac
{
2l |s|
}
{
(1+l^{2})a+(1-l^{2}) |s|
}
-\frac{k|s|}{a}.
\end{equation}
Then we conclude 
\begin{equation*}\label{value of l}
l
\leq
\frac
{2ka +(1-k^{2})|b|}
{(1+k^{2})a + (1-k^{2})\sqrt{a^{2}+b^{2}}}.
\end{equation*}
which is suitable for our purpose.

(ii)
In order to eliminate the additional assumption that $f(z)/z \neq 0$ in $\D$, we need a sort of stability of the condition \eqref{maineq};

\begin{lem}\label{approx}
If $f \in \A$ satisfies the assumption of Theorem \ref{main}, then so does $\dstyle f_{r}(z) = \frac1r f(rz),\,r \in (0,1)$.
\end{lem}

\pf
It follows from the assumption that $aH_{s}(rz)$ is contained in the disk
$$
\Delta = 
\left\{ w \in \C :
\left|
w -\frac{cr^{2}|z|^{2}+s}{1-r^{2}|z|^{2}}
\right|
\leq
\frac{M}{1-r^{2}|z|^{2}}
\right\}.
$$
We want to deduce that $aH_{s}(rz)$ lies in the disk
$$
\Delta' = 
\left\{ w \in \C :
\left|
w -\frac{c|z|^{2}+s}{1-|z|^{2}}
\right|
\leq
\frac{M}{1-|z|^{2}}
\right\}.
$$
Therefore it is enough to see that $\Delta \subset \Delta'$, that is,
\begin{equation}\label{circle3}
\left|
\frac{c|z|^{2}+s}{1-|z|^{2}}  -  \frac{cr^{2}|z|^{2}+s}{1-r^{2}|z|^{2}}
\right|
\leq
\frac{M}{1-|z|^{2}}  -  \frac{M}{1-r^{2}|z|^{2}}.
\end{equation}
In view of the identity
$$
\frac{|z|^{2}}{1-|z|^{2}} - \frac{r^{2}|z|^{2}}{1-|z|^{2}}
=
\frac{1}{1-|z|^{2}} - \frac{1}{1-r^{2}|z|^{2}},
$$
the inequality \eqref{circle3} is equivalent to \eqref{inequality}.
\epf

Now we shall show that the condition $f(z)/z \neq 0$ in $\D$ follows from the assumption of Theorem \ref{main}.
Suppose, to the contrary, that $f(z_{0}) =0$ for some $0 < |z_{0}| < 1$.
We may assume that $f(z) \neq 0$ for $0 < |z| < |z_{0}|$.
Then by Lemma \ref{approx} we can apply Theorem \ref{main} to the function $f_{r_{0}}(z) = f(r_{0}z)/r_{0},\, r_{0}=|z_{0}|$ to conclude that $f_{r_{0}}$ has a quasiconformal extension to $\C$.
In particular, $f_{r_{0}}$ is injective on $\overline{\D}$.
This, however, contradicts the relation $f_{r_{0}}(z_{0}/r_{0}) = f_{r_{0}}(0) = 0$.
\epf


\begin{rem}
We can replace $|s|$ in \eqref{chain2} to any positive real value and continue our argument.
However, it will be found that $|s|$ gives the smallest $l$ by calculations.
\end{rem}

\begin{rem}
We have $l \geq k$, where $l=k$ if and only if $b = 0$.
Indeed, let $l=l(k)$. 
Then we have $l'(k) > 0$ and $l''(k) \leq 0$ which imply $l \geq k$.
If we suppose $l=k \neq 0$, then the right-hand side of \eqref{Inequality01} is greater than or equal to $0$ only if $b=0$.
In the case $l=k=0$ we also have $b=0$ by \eqref{Inequality01}.
It easily follows from \eqref{mainl} that $l=k$ if $b=0$.
\end{rem}


\

\section{Proof of Theorem \ref{maincor}}
It is easy to see from \eqref{bazilevic} that $f$ is a Bazilevi\v c function of type ($\alpha, \beta$) under our assumption since $M$ is always less than or equal to $(\alpha^{2}+\beta^{2})/\alpha$.

Let us now prove quasiconformal extensibility of $f$.
Setting $1/s  =\alpha+i \beta$ which implies $a = \textup{Re}s =\alpha/(\alpha^{2}+\beta^{2})$ and $b=\textup{
Im}s =-\beta/(\alpha^{2}+\beta^{2})$, \eqref{maincorineq} turns to
$$
\left| 
1+ \frac{zf''(z)}{f'(z)}+\left(\frac1s-1\right)\frac{zf'(z)}{f(z)}-\frac1a
\right|
\leq 
\left\{
\begin{array}{ll}
 k , &  0 < a < 1, \\[5pt]
k/a , &  1 \leq a.
\end{array} 
\right.
$$
Therefore, Theorem \ref{maincor} follows from Theorem \ref{main} with $c=-s$.
\epf


\

\section{Proof of Theorem \ref{mainD*}}
First let $s \neq 1$.
In that case we may assume $g(\zeta) \neq 0$ for all $\zeta \in \D^{*}$ because of a similar discussion of the proof of Theorem \ref{main};

\begin{lem}\label{approxD*}
Let $\dstyle g(\zeta)= \zeta +\frac{d}{\zeta} + \cdots$ be analytic in $\D^{*}$.
If $g$ satisfies the same assumption of Theorem \ref{mainD*}, then so does $\dstyle g_{R}(\zeta) = \frac{1}{R}f(R\zeta)$,\,$R>1$.
\end{lem}
\pf
We need to prove
$$
\left|
\frac{ib}{|\zeta|^{2}-1} - a G_{s}(R\zeta)
\right|
\leq
\frac{ak|s|-|b|(a-1)}{|\zeta|^{2}-1}
$$ 
by using
$$
\left|
\frac{ib}{R^{2}|\zeta|^{2}-1} - a G_{s}(R\zeta)
\right|
\leq
\frac{ak|s|-|b|(a-1)}{R^{2}|\zeta|^{2}-1},
$$
where
$$
G_{s}(\zeta) = (1-s) \left( \frac{\zeta g'(\zeta)}{g(\zeta)}-1\right) + s\frac{\zeta g''(\zeta)}{g'(\zeta)}.
$$
In a similar way to the proof of Lemma \ref{approx}, it suffices to see that
$$
\left|
\frac{ib}{|\zeta|^{2}-1} - \frac{ib}{R^{2}|\zeta|^{2}-1}
\right|
\leq
\frac{ak|s|-|b|(a-1)}{|\zeta|^{2}-1} - \frac{ak|s|-|b|(a-1)}{R^{2}|\zeta|^{2}-1}.
$$
This is equivalent to $|b| \leq k|s|$.
\epf

Then we let
$$
f(1/\zeta,t) = 
\frac{1}{g(e^{st}\zeta)}
\left\{
1-(1-e^{-2t})  e^{st}\zeta
\frac{g'(e^{st}\zeta)}{g(e^{st}\zeta)}
\right\}^{-s}
$$
and
$$
F(1/\zeta,t) = f(1/\zeta,t/|s|).
$$
Since
$$
h(1/\zeta,t)
= 
\frac{\dot{F}(1/\zeta,t)}{(1/\zeta) F'(1/\zeta,t)}
=
\frac{s}{|s|}\cdot
\frac{1+P(e^{st/|s|}\zeta,t/|s|)}{1-P(e^{st/|s|}\zeta,t/|s|)}
$$
where
$$
P(\zeta,t) = (e^{2t/|s|}-1) G_{s}(\zeta),
$$
it is sufficient to see that
\begin{equation}\label{diskBD*}
|a P(e^{st/|s|}\zeta,t/|s|) + ib| <  k|s|
\end{equation}
for all $\zeta \in \D^{*}$ and $t \in [0,\infty)$ under the assumption of the theorem. 
By triangle inequality we have
$$
|aP+ib|
\leq
\left|
\frac{1-e^{2t/|s|}}{1-e^{2at/|s|}}
\left(ib+(1-e^{2at/|s|})aG_{s}(e^{st/|s|}\zeta)\right)
\right|
+
\left|
ib\left( 1-\frac{1-e^{2t/|s|}}{1-e^{2at/|s|}}\right)
\right|
$$
for $\zeta \in \D^{*}$ and $t \in [0,\infty)$.
Following the lines of the proof of Lemma \ref{lemdiskB}, one can obtain that \eqref{maineqD*} implies \eqref{diskBD*}.
Therefore, a similar argument of the proof of Theorem \ref{main} implies our assertion.
The case $s=1$ follows from a theorem of Becker \cite{Becker:1972}.
\epf


\

\section{Proof of Corollary \ref{corD*} and \ref{corD*D}}
\noindent{\bf\textit{Proof of Corollary \ref{corD*}.}}
Let $R>1$ be an arbitrary but fixed number.
We would like to show that $g_{R}(\zeta) = g(R\zeta)/R$ can be extended to a $k$-quasiconformal mapping of $\CC-\{0\}$.
Since $g(\zeta) \neq 0$ in $\zeta \in \D^{*}$ from the assumption, there exists a certain constant $A$ such that
$$
(|\zeta|^{2}-1)\left|1-\frac{\zeta g_{R}'(\zeta)}{g_{R}(\zeta)}\right|   \leq A < \infty
$$
for all $\zeta \in \overline{\D^{*}}$.
We also have
$$
\left|
1-
\frac{\zeta g_{R}'(\zeta)}{g_{R}(\zeta)}+
\frac{\zeta g_{R}''(\zeta)}{g_{R}'(\zeta)}
\right|  \leq
\frac{k}{|\zeta R|^{2}-1}
$$
for $\zeta \in \D^{*}$.
Thus we obtain with $s=R^{2}A/k(R^{2}-1)$
$$
(|\zeta|^{2}-1)
\left|
\frac{1}{s}
\left(
1-\frac{\zeta g_{R}'(\zeta)}{g_{R}(\zeta)}
\right)
-1-
\frac{\zeta g_{R}''(\zeta)}{g_{R}'(\zeta)} +
\frac{\zeta g_{R}'(\zeta)}{g_{R}(\zeta)}
\right|
\leq
\frac{A}{s}+
k\frac{|\zeta|^{2}-1}{|\zeta R|^{2}-1}
\leq k
$$
which implies quasiconformal extensibility of $g_{R}$ by Theorem \ref{mainD*}.
A limiting procedure proves Corollary \ref{corD*}. 
\epf

\noindent{\bf\textit{Proof of Corollary \ref{corD*D}.}}
Note that the function
$
1 + (zf''(z)/f'(z)) - (zf'(z)/f(z))
$
is analytic in $\D$ and has a zero of order 2 at the origin by the condition $f''(0)=0$.
Thus, we obtain from the assumption that
$$
\frac{1}{|z|^{2}}
(1-|z|^{2})
\left|
1 + \frac{zf''(z)}{f'(z)} - \frac{zf'(z)}{f(z)}
\right|
\leq k
$$
by the maximum modulus principle.
Let $g(\zeta)$ be a function defined by
$$
g(\zeta) = \frac{1}{f(z)}
$$
where $\zeta = 1/z$. 
Then $g$ is analytic in $\D^{*}$ and has the form $g(\zeta) = \zeta + d/\zeta + \cdots$.
From the relations
$$
\frac{zf'(z)}{f(z)} = \frac{\zeta g'(\zeta)}{g(\zeta)}
$$
and
$$
1 + \frac{zf''(z)}{f'(z)} = -1-\frac{\zeta g''(\zeta)}{g'(\zeta)} +2\frac{\zeta g'(\zeta)}{g(\zeta)},
$$
we can deduce our assertion by applying Corollary \ref{corD*}.
\epf


\

\no
\textbf{Acknowledgement.} 
The author expresses his deep gratitude to Professor Toshiyuki Sugawa for many helpful discussions and his continuous encouragement during this work.

\bibliographystyle{amsplain}
\bibliography{bibdata}
\end{document}